\documentclass[a4paper]{amsart}

\usepackage{amsmath,amsthm,amsfonts,mathrsfs}

\newtheorem{theorem}{Theorem}
\newtheorem{lemma}[theorem]{Lemma}
\newtheorem{corollary}[theorem]{Corollary}

\theoremstyle{remark}

\DeclareMathOperator{\probability}{\mathsf{P}}
\DeclareMathOperator{\expectation}{\mathsf{E}}
\DeclareMathOperator{\variance}{\mathsf{Var}}

\allowdisplaybreaks[3]

\title{A Note on Random Coverings of Tori}

\author{Tomas Persson} \address{Tomas Persson\\ Centre for
  Mathematical Sciences\\ Lund University\\ Box 118\\ 22100
  Lund\\ Sweden} \email{tomasp@maths.lth.se}

\begin{document}

\begin{abstract}
  This note provides a generalisation of a recent result by
  J\"arven\-p\"a\"a, J\"ar\-ven\-p\"a\"a, Koivusalo, Li, and Suomala,
  (to appear), on the dimension of limsup-sets of random coverings of
  tori. The result in this note is stronger in the sense that it
  provides also a large intersection property of the limsup-sets, the
  assumptions are weaker, and it implies the result of
  J\"arven\-p\"a\"a, J\"ar\-ven\-p\"a\"a, Koivusalo, Li, and Suomala
  as a special case. The proof is based on a recent result by Persson
  and Reeve from 2013.
\end{abstract}


\maketitle

\section{Introduction}

Let $d$ be a natural number. We consider the $d$-dimensional torus
$\mathbb{T}^d$, and a sequence of open sets $U_i \subset
\mathbb{T}^d$. The random vectors $v_i$ are independent and uniformly
distributed on the torus $\mathbb{T}^d$, and are used to translate the
sets $U_i$, hence producing a sequence $V_i (v_i)$ of random sets
defined by $V_i (v_i) = U_i + v_i$. We are interested in the typical
behaviour of the limsup-set
\[
E (v) = \limsup_{i \to \infty} V_i (v_i),
\]
that is, the set of points on the torus that are covered by
infinitely many sets $V_i (v_i)$.

Limsup-sets often possess a large intersection property, see Falconer
\cite{falconer}. This means that the set belongs, for some $0 < s \leq
d$, to the class $\mathscr{G}^s (\mathbb{T}^d)$, where $\mathscr{G}^s
(\mathbb{T}^d)$ is the largest collection of $G_\delta$ subsets of
$\mathbb{T}^d$ with the property that any countable intersection of
such sets has Hausdorff dimension at least $s$. For instance, we have
$\mathscr{G}^s (\mathbb{T}^d) \subset \mathscr{G}^t (\mathbb{T}^d)$
provided $t < s$, and if $A \in \mathscr{G}^t (\mathbb{T}^d)$ for all
$t < s$, then $A \in \mathscr{G}^s (\mathbb{T}^d)$. For more
properties of these classes, relevant in this paper, we refer the
reader to the paper \cite{perssonreeve}. In this note, we shall be
concerned with the large intersection properties of typical $E (v)$.

Let $\lambda$ denote the $d$-dimensional Lebesgue measure on
$\mathbb{T}^d$. For $0 < s < d$ and a set $A \subset \mathbb{T}^d$, we
define the $s$-energy of $A$ as
\[
I_s (A) = \iint_{A \times A} |x - y|^{-s}\, \mathrm{d}x \mathrm{d}y,
\]
where $|x - y|$ denotes the distance between the points $x$ and $y$.

The aim of this note is to give a short proof of the following
theorem. For the background of this and other similar results on
random coverings of tori, we refer the reader to
\cite{jarvenpaa_etal}.

\begin{theorem} \label{the:randomcovering}
  The set $E(v)$ is almost surely in the class $\mathscr{G}^s
  (\mathbb{T}^d)$, where $s$ is defined by
  \[
  s = \inf \{\, t : \sum_{i=1}^\infty \frac{\lambda (U_i)^2}{I_t
    (U_i)} < \infty \text{ or } t = d \,\}.
  \]
\end{theorem}

In the paper \cite{jarvenpaa_etal}, J\"arvenp\"a\"a, J\"arvenp\"a\"a,
Koivusalo, Li, and Suomala proved a similar result. They imposed more
restrictive assumptions on the sets, and they only proved the
dimension result, not the large intersection property. It is not
immediately clear if the result in \cite{jarvenpaa_etal} provides the
same dimension result that Theorem~\ref{the:randomcovering} does,
under the extra conditions imposed in \cite{jarvenpaa_etal}. However,
we shall study below two corollaries of
Theorem~\ref{the:randomcovering}. The second corollary will show that
the result of J\"arvenp\"a\"a, J\"arvenp\"a\"a, Koivusalo, Li, and
Suomala is a special case of Theorem~\ref{the:randomcovering}. Hence,
this note generalises the paper \cite{jarvenpaa_etal}, providing a
stronger result under weaker assumptions. Moreover, the proof is much
shorter.

To derive corollaries of Theorem~\ref{the:randomcovering}, we will
estimate the $t$-energies $I_t (U_i)$. For the first corollary, we do
this as follows. If $B_i = B_i (0,r_i)$ is a ball with $\lambda (U_i)
= \lambda (B_i)$, then we may estimate that
\[
I_t (U_i) \leq I_t (B_i) = C_t r_i^{2d - t} = C_t' \lambda
(B_i)^{2-t/d} = C_t' \lambda (U_i)^{2 - t/d}.
\]
where $C_t$ and $C_t'$ are constants. Hence,
\[
\sum_{i=1}^\infty \lambda (U_i)^{t/d} = \infty \qquad \Rightarrow
\qquad \sum_{i=1}^\infty \frac{\lambda (U_i)}{I_t (U_i)} = \infty,
\]
and we get the following corollary to
Theorem~\ref{the:randomcovering}.

\begin{corollary} \label{cor:randomcovering}
  The set $E(v)$ is almost surely in the class $\mathscr{G}^s
  (\mathbb{T}^d)$, where $s$ is defined by
  \[
  s = \inf \{\, t : \sum_{i=1}^\infty \lambda (U_i)^{t/d} < \infty
  \text{ or } t = d \,\}.
  \]
\end{corollary}

Corollary~\ref{cor:randomcovering} does not always provide the optimal
result, whereas the result in \cite{jarvenpaa_etal} does in the case
considered there. To clarify the differences, let us study an
example. Let $d = 2$. Suppose $1 < \alpha < \beta$, and that $U_i$ is
a rectangle with side lengths about $1/i^\alpha$ and $1/i^\beta$. Then
$\lambda (U_i) = 1/i^{\alpha + \beta}$, and
Corollary~\ref{cor:randomcovering} implies that almost surely $E (v)$
is in the class $\mathscr{G}^{2/(\alpha + \beta)}
(\mathbb{T}^d)$. However, by J\"arvenp\"a\"a, J\"arvenp\"a\"a,
Koivusalo, Li, and Suomala, the dimension is almost surely
$1/\alpha$. Since $\frac{2}{\alpha + \beta} < \frac{1}{\alpha}$, this
shows that Corollary~\ref{cor:randomcovering} does not give the
optimal result, (at least not when it comes to dimension).

Note however that in the case $d=1$,
Corollary~\ref{cor:randomcovering} gives the optimal result. In this
case it was proved by Durand when $U_i$ are intervals \cite{durand}.

The reason that Corollary~\ref{cor:randomcovering} is not optimal is
that if the sets $U_i$ are not sufficiently similar to balls, then it
is to rough an estimate to estimate $U_i$ by the ball $B_i$, as was
done above. If $U_i$ is comparable to a $d$-dimensional rectangle, as
in \cite{jarvenpaa_etal}, then one would do better estimating $U_i$ by
such a rectangle. We shall do so in what follows.

Suppose $Q$ is a $d$-dimensional cube, and for each $i$ we have that
$R_i = L_i (Q) \subset U_i$, where $L_i$ is an affine transformation
with singular values
\[
0 < \alpha_d (L_i) \leq \cdots \leq \alpha_1 (L_i) < 1.
\]
We define as in \cite{jarvenpaa_etal}, the singular value function
\[
\Phi^s (L_i) = \alpha_1 (L_i) \alpha_2 (L_i) \cdots \alpha_{m-1} (L_i)
\alpha_m^{s-m+1} (L_i),
\]
where $m$ is such that $m - 1 < s \leq m$.

One can easily show that in this case, there is a constant $K$ such
that
\[
I_s (R_i) \leq K \frac{\lambda (R_i)^2}{\Phi^s (L_i)}.
\]
Hence we get that
\[
\sum_{i=1}^\infty \frac{\lambda (R_i)^2}{I_s (R_i)} \geq K^{-1}
\sum_{i=1}^\infty \Phi^s (L_i).
\]
This gives us the following corollary of
Theorem~\ref{the:randomcovering}. It is essentially the result in
\cite{jarvenpaa_etal}, but it is stronger since it also gives the
large intersection property, and imposes somewhat less restrictive
assumptions.

\begin{corollary}
If $R_i \subset U_i$ as above, then the set $E (v)$ is almost surely
in the class $\mathscr{G}^s (\mathbb{T}^d)$, where $s$ is defined by
\[
s = \inf \{\, t: \sum_{i=1}^\infty \Phi^s (L_i) < \infty \text{ or } t
= d \, \}.
\]
\end{corollary}

\section{Proof of Theorem \ref{the:randomcovering}}


The proof is based on the following lemma from \cite{perssonreeve},
that gives us a method to determine if a limsup-set belongs to the
class $\mathscr{G}^s (\mathbb{T}^d)$. The theorem is only stated and
proved for $d=1$ in \cite{perssonreeve}, but it holds for any $d$, and
only minor changes in the proof are required to make it work for $d >
1$. Also, the statement in \cite{perssonreeve} is for $[0,1]$ instead
of $\mathbb{T}^1$, but this difference is not substantial.

\begin{lemma} \label{the:frostman}
  Let $E_k$ be open subsets of\/ $\mathbb{T}^d$, and $\mu_k$ Borel
  probability measures, with support in the closure of $E_k$, that
  converge weakly to a measure $\mu$ with density $h$ in $L^2$. Assume
  that $\mu (I) > 0$ for all cubes $I \subset [0, 1)^d$ with non-empty
    interior, and assume that for each $\varepsilon > 0$, there is a
    constant $C_\varepsilon$, such that
  \begin{equation} \label{eq:est_on_limitmeasure}
    |I|^{1+\varepsilon} \lVert h \chi_I \rVert_2^2 \leq C_\varepsilon
    \lVert h \chi_I \rVert_1^2
  \end{equation}
  holds for any cube $I \subset \mathbb{T}^d$. If there is a constant
  $C$ such that
  \begin{equation} \label{eq:bounded_energy}
    \iint |x - y|^{-s} \, \mathrm{d} \mu_k (x) \mathrm{d} \mu_k (y)
    \leq C
  \end{equation}
  holds for all $k$, then $\limsup E_k$ is in the class
  $\mathscr{G}^s (\mathbb{T}^d)$.
\end{lemma}

In our application of Lemma~\ref{the:frostman}, the limit measure
$\mu$ will be the Lebesgue measure, and therefore the assumption
\eqref{eq:est_on_limitmeasure} will be automatically fulfilled. Note
also that the proof of Lemma~\ref{the:frostman} can be significantly
simplified in this case.

Let $E_k (v) = \bigcup_{i = m_k}^k V_i (v_i)$, where $m_k < k$ is a
sequence increasing to infinity. We then have $\limsup E_k (v) = E (v)
= \limsup V_i (v_i)$.  Define $\mu_k = \sum_{i=m_k}^k c_{i,k}
\lambda|_{V_i (v_i)}$, where $c_{i,k}$ are constants that will be
specified later, but are such that $\mu_k$ are probability
measures. In particular, $\sum_{i=m_k}^k \sum_{j=m_k}^k c_{i,k}
c_{j,k} \lambda (U_i) \lambda (U_j) \leq 1$.

Let $s = \inf \{\, t : \sum_i \lambda (U_i)^2 / I_t (U_i) < \infty
\,\}$, and pick $t$ with $t < s$ and $t < d$. We need to prove that
with probability 1, we have $E (v) \in \mathscr{G}^t (\mathbb{T}^d)$.

If $i \neq j$ we have, since $v_i$ and $v_j$ are independent and
uniformly distributed, that
\begin{equation} \label{eq:energy_independent}
  \expectation \biggl( \iint_{V_i (v_i) \times V_j (v_j)} |x - y|^{-t}
  \, \mathrm{d}x \mathrm{d}y \biggr) \leq C \lambda (U_i) \lambda
  (U_j),
\end{equation}
where $C$ is a constant that only depends on $t$ and
$d$. ($\expectation$ denotes expectation.)

However, if $i = j$, then $v_i$ and $v_j$ are not at all
independent. We then have
\begin{equation} \label{eq:energy_dependent}
  \expectation \biggl( \iint_{V_i (v_i) \times V_i (v_i)} |x - y|^{-t}
  \, \mathrm{d}x \mathrm{d}y \biggr) = I_t (U_i).
\end{equation}

We now use the estimates \eqref{eq:energy_independent} and
\eqref{eq:energy_dependent}, to get that
\begin{align*}
  \expectation \biggl( \iint |x - y|^{-t} \, \mathrm{d}\mu_k(x)
  \mathrm{d}\mu_k(y) \biggr) &\leq \sum_{i \neq j} C c_{i,k} c_{j,k}
  \lambda(U_i) \lambda(U_j) + \sum_{i=m_k}^k c_{i,k}^2 I_t (U_i) \\ &
  \leq C + \sum_{i=m_k}^k c_{i,k}^2 I_t (U_i).
\end{align*}

Hence, to be able to apply Lemma~\ref{the:frostman}, we need to
choose the numbers $c_{i,k}$ and a number $L$ so that
\[
\sum_{i=m_k}^k c_{i,k} \lambda (U_i) = 1 \quad \text{and} \quad
\sum_{i=m_k}^k c_{i,k}^2 I_t (U_i) \leq L \quad \text{for all } k,
\]
and $\mu_k$ almost surely converges weakly to the $d$-dimensional
Lebesgue measure, at least along a sub sequence.

Let $c_{i,k} = c_k \lambda (U_i) / I_t (U_i)$, with $c_k = \bigl(
\sum_{i=m_k}^k \lambda (U_i)^2 / I_t (U_i) \bigr)^{-1}$. Then
\[
\sum_{i=m_k}^k c_{i,k} \lambda (U_i) = c_k \sum_{i=m_k}^k
\frac{\lambda (U_i)^2}{I_t (U_i)} = 1,
\]
and
\[
\sum_{i=m_k}^k c_{i,k}^2 I_t (U_i) = c_k^2 \sum_{i=m_k}^k
\frac{\lambda (U_i)^2}{I_t (U_i)} = c_k.
\]
Because of the choice of $t$, we have that $c_k$ converges to $0$ as
$k$ grows, provided $m_k$ grows sufficiently slow. We may therefore
choose $L = \sup c_k < \infty$.

Finally, we observe that the fact that $c_k \to 0$ as $k \to \infty$,
implies that there is a sequence $n_k$ such that $\mu_{n_k}$ almost
surely converges weakly to the Lebesgue measure. To see this, pick any
continuous function $\phi \colon \mathbb{T}^d \to \mathbb{R}$, and
define the random variables $X_i = \int_{V_i (v_i)} \phi \,
\mathrm{d}\lambda$ and $S_k = \sum_{i=m_k}^k c_{i,k} X_i$.  Then
\[
\expectation S_k = \lambda (\phi) := \int \phi \, \mathrm{d}\lambda,
\qquad \text{and} \qquad \variance S_k = \sum_{i = m_k}^k c_{i,k}^2
\variance X_i.
\]
Since
\[
\variance X_i \leq \bigl((\sup \phi)^2 - \lambda (\phi)^2 \bigr)
\lambda (U_i)^2 = C_\phi \lambda(U_i)^2,
\]
and $\lambda(U_i)^2 / I_t (U_i) \leq 1$, we derive that $\variance S_k
\leq C_\phi c_k \to 0$. We choose any sub sequence $n_k$ of the
natural numbers with $n_k \to \infty$ and $\sum_{k=1}^\infty c_{n_k} <
\infty$. Then for any $\varepsilon > 0$
\[
\sum_{k = m}^\infty \probability (|S_{n_k} - \lambda (\phi)| >
\varepsilon) \leq \sum_{k=m}^\infty \frac{\variance
  S_{n_k}}{\varepsilon^2} \to 0, \quad m \to \infty.
\]
This implies that $S_{n_k} \to \int \phi \, \mathrm{d} \lambda$ almost
surely, and therefore we have that almost surely $\mu_{n_k}$ converges
weakly to $\lambda$.

Lemma~\ref{the:frostman} now finishes the proof.

\section*{Acknowledgement}

The author is grateful to Esa J\"arvenp\"a\"a for pointing out a
mistake.

\end{document}